\def\pf{{\bf Proof:}}

\def\beg{\begin{equation}}
\def\en{\end{equation}}
\def\beq{\begin{eqnarray}}
\def\enq{\end{eqnarray}}
\def\ex{exceptional\ }
\def\vv{vector\ }

\def\bbbr{{\rm I\!R}} 
\def\bbbn{{\rm I\!N}} 

\def\bbbp{{\rm I\!P}}

\def\bbbc{{\mathchoice {\setbox0=\hbox{$\displaystyle\rm C$}\hbox{\hbox
to0pt{\kern0.4\wd0\vrule height0.9\ht0\hss}\box0}}
{\setbox0=\hbox{$\textstyle\rm C$}\hbox{\hbox
to0pt{\kern0.4\wd0\vrule height0.9\ht0\hss}\box0}}
{\setbox0=\hbox{$\scriptstyle\rm C$}\hbox{\hbox
to0pt{\kern0.4\wd0\vrule height0.9\ht0\hss}\box0}}
{\setbox0=\hbox{$\scriptscriptstyle\rm C$}\hbox{\hbox
to0pt{\kern0.4\wd0\vrule height0.9\ht0\hss}\box0}}}}

\def\bbbg{{\mathchoice {\setbox0=\hbox{$\displaystyle\rm G$}\hbox{\hbox
to0pt{\kern0.4\wd0\vrule height0.9\ht0\hss}\box0}}
{\setbox0=\hbox{$\textstyle\rm G$}\hbox{\hbox
to0pt{\kern0.4\wd0\vrule height0.9\ht0\hss}\box0}}
{\setbox0=\hbox{$\scriptstyle\rm G$}\hbox{\hbox
to0pt{\kern0.4\wd0\vrule height0.9\ht0\hss}\box0}}
{\setbox0=\hbox{$\scriptscriptstyle\rm G$}\hbox{\hbox
to0pt{\kern0.4\wd0\vrule height0.9\ht0\hss}\box0}}}}

\def\bbbq{{\mathchoice {\setbox0=\hbox{$\displaystyle\rm Q$}\hbox{\raise
0.15\ht0\hbox to0pt{\kern0.4\wd0\vrule height0.8\ht0\hss}\box0}}
{\setbox0=\hbox{$\textstyle\rm Q$}\hbox{\raise
0.15\ht0\hbox to0pt{\kern0.4\wd0\vrule height0.8\ht0\hss}\box0}}
{\setbox0=\hbox{$\scriptstyle\rm Q$}\hbox{\raise
0.15\ht0\hbox to0pt{\kern0.4\wd0\vrule height0.7\ht0\hss}\box0}}
{\setbox0=\hbox{$\scriptscriptstyle\rm Q$}\hbox{\raise
0.15\ht0\hbox to0pt{\kern0.4\wd0\vrule height0.7\ht0\hss}\box0}}}}
\def\bbbt{{\mathchoice {\setbox0=\hbox{$\displaystyle\rm
T$}\hbox{\hbox to0pt{\kern0.3\wd0\vrule height0.9\ht0\hss}\box0}}
{\setbox0=\hbox{$\textstyle\rm T$}\hbox{\hbox
to0pt{\kern0.3\wd0\vrule height0.9\ht0\hss}\box0}}
{\setbox0=\hbox{$\scriptstyle\rm T$}\hbox{\hbox
to0pt{\kern0.3\wd0\vrule height0.9\ht0\hss}\box0}}
{\setbox0=\hbox{$\scriptscriptstyle\rm T$}\hbox{\hbox
to0pt{\kern0.3\wd0\vrule height0.9\ht0\hss}\box0}}}}
\def\bbbs{{\mathchoice
{\setbox0=\hbox{$\displaystyle     \rm S$}\hbox{\raise0.5\ht0\hbox
to0pt{\kern0.35\wd0\vrule height0.45\ht0\hss}\hbox
to0pt{\kern0.55\wd0\vrule height0.5\ht0\hss}\box0}}
{\setbox0=\hbox{$\textstyle        \rm S$}\hbox{\raise0.5\ht0\hbox
to0pt{\kern0.35\wd0\vrule height0.45\ht0\hss}\hbox
to0pt{\kern0.55\wd0\vrule height0.5\ht0\hss}\box0}}
{\setbox0=\hbox{$\scriptstyle      \rm S$}\hbox{\raise0.5\ht0\hbox
to0pt{\kern0.35\wd0\vrule height0.45\ht0\hss}\raise0.05\ht0\hbox
to0pt{\kern0.5\wd0\vrule height0.45\ht0\hss}\box0}}
{\setbox0=\hbox{$\scriptscriptstyle\rm S$}\hbox{\raise0.5\ht0\hbox
to0pt{\kern0.4\wd0\vrule height0.45\ht0\hss}\raise0.05\ht0\hbox
to0pt{\kern0.55\wd0\vrule height0.45\ht0\hss}\box0}}}}
\def\bbbz{{\mathchoice {\hbox{$\sf\textstyle Z\kern-0.4em Z$}}
{\hbox{$\sf\textstyle Z\kern-0.4em Z$}}
{\hbox{$\sf\scriptstyle Z\kern-0.3em Z$}}
{\hbox{$\sf\scriptscriptstyle Z\kern-0.2em Z$}}}}

\def\gaaa{\ifmmode
              {{\mbox{\deu a}}}%
          \else${{\mbox{\deu a\ }}}$%
           \fi}
\def\gbbb{\ifmmode
              {{\mbox{\deu b}}}%
          \else${{\mbox{\deu b\ }}}$%
           \fi}
\def\ccc{\ifmmode
        {\bbbc}%
          \else${\bbbc\ }$%
          \fi}
\def\fff{\ifmmode
        {{\mbox{\bf F}}}%
          \else${{\mbox{\bf F\ }}}$%
          \fi}
\def\fxfx{\ifmmode
              {{\bbbf}_p[x] }%
          \else${{\bbbf}_p[x]\ }$%
          \fi}
\def\fffp{\ifmmode
              {{\mbox{\bf F}}_p }%
          \else${{{\mbox{\bf F}}}_p\ }$%
          \fi}
\def\fffq{\ifmmode
              {{\bbbf}_q }%
          \else${{\bbbf}_q\ }$%
          \fi}
\def\gmmm{\ifmmode
              {{\mbox{\deu m}}}%
          \else${{\mbox{\deu m\ }}}$%
           \fi}
\def\nnn{\ifmmode
        {\bbbn}%
          \else${\bbbn\ }$%
          \fi}

\def\nnno{\ifmmode
        {{\bbbn}_0}%
          \else${{\bbbn}_0\ }$%
          \fi}

\def\ooo{\ifmmode
              {{\mbox{\deu o}}}%
          \else${{\mbox{\deu o\ }}}$%
           \fi}
\def\gppp{\ifmmode
              {{\mbox{\deu p}}}%
          \else${{\mbox{\deu p\ }}}$%
           \fi}

\def\ppp{\ifmmode
              {{\mbox{\deu p}}}%
          \else${{\mbox{\deu p\ }}}$%
           \fi}
\def\cppp{\ifmmode
              {\cal P}%
          \else${\cal P\ }$%
           \fi}

\def\qqq{\ifmmode
              {\bbbq}%
          \else${\bbbq\ }$%
           \fi}
\def\zzz{\ifmmode
              {\bbbz}%
          \else${\bbbz\ }$%
          \fi}
\def\grrr{\ifmmode
              {{\mbox{\deu r}}}%
          \else${{\mbox{\deu r\ }}}$%
           \fi}
\def\rrr{\ifmmode
              {\bbbr}%
          \else${\bbbr\ }$%
          \fi}
\def\st{\ifmmode
              {***}%
         \else${***}$%
          \fi}

\newcommand{\deu}{\bf}

\newcommand{\text}{\mbox}

\documentstyle[12pt]{article}
\begin{document}
\title{ Exceptional vector bundle on Enriques surfaces}
\author{Severinas Zube}
\date{1994 October}
\maketitle
\begin{abstract}
The main purpose in this paper is to study exceptional vector bundles
on  Enriques surfaces.
\end{abstract}

\subsection{ Introduction}
The purpose of this note is to study \ex \vv bundles on Enriques surfaces.
Exceptional bundle E on a surface with irregularity $q=h^1(O_S)$ and geometric
genus
 $p_g=h^2(O_S)$ is the bundle with the following properties:\\
$Ext^0(E,E)=\bbbc
,Ext^1(E,E)=q, Ext^2(E,E)=p_g$. On an Enriques surface Kim and Naie in
[Ki1],[Ki2],[N] have been studied extremal bundles which are very similar
to \ex\ bundles. Extremal bundle on Enriques surface by definition is a
simple $Ext^0(E,E)=\bbbc$, rigid $Ext^1(E,E)=\bbbc$ with the following
condition:$Ext^2(E,E)=\bbbc$. From the Riemann-Roch theorem easily follows
that any \ex bundle has odd rank and any extremal even rank. In [Ki2],
Kim characterized  extremal bundles on Enriques surfaces. They exist only
on nodal surfaces and satisfy $c_1^2=4n-2,c_2=n$ for $n>4$. They have also
geometrical meaning. It turns out that the existence of extremal bundles
is closely related to the embedding a general Enriques surface in the
Grassmannian G(2,n+1).\\
     The main result in this paper is to give the necessary and sufficient
conditions for the existence \ex bundles on Enriques surfaces. The
statement  is similar to the theorem 4 in [Ki2]. The proof fill  a
gap in the proof of this theorem 4. Also I give constructions of them
by using some
constructions of Enriques surfaces and by using modular operations
(reflections) which is
described in the last section. I think the reflection is very useful to
construct and study moduli of sheaves on Enriques surfaces.

\subsection{Enriques surfaces}
A smooth irreducible surface $S$, such that $h^1 (O_S )=h^2 (O_S )=0 $ and
$2K_S \sim O_S $, is called a Enriques surface.
Recall that a divisor $D$ on a smooth surface $X$ is said to be $nef$
 if $DC \geq 0$ for every curve $C$ on $X$ . The following useful properties
 will be  used throughout, sometimes without explicit mention:

(A)([C,D] Corollary 3.1.3) If $D$ is a $nef$ divisor  and $D^2 > 0 $ , then
$H^1 (O_S (-D)) = 0$ and $\chi(O(D)) -1 =dim|D| = \frac{D^2 }{2}$. \label{222}

(B)([C,D] Proposition 3.1.4) If $ \mid D \mid $ has no fixed components,   then
one of the following holds: \nonumber \\
(i) $D^2 > 0$ and there exist an irreducible curve $C$ in $ \mid D \mid $.
\nonumber \\
(ii) $D^2 =0$ and there exist a genus 1 pencil $ \mid P \mid $ such that $D
\sim kP$ for some $k \geq 1 $. \label{333}

(C)([C,D]Chapter 4, appendix , corollary 1. and corollary 2.) If $D^2 \geq 6$
and $D$ is $nef$ then $D$ is ample , $2D$ is generated by its global sections,
$3D$ is very ample.

The Enriques surface $S$ is called nodal (resp. unnodal) if there are
(resp. not) a smooth (-2)-curve
contained in $S$. A general Enriques surface is an unnodal.
Every Enriques surface admits an elliptic fibration over $\bbbp ^1$
with exactly two multiple fibers F,  F' and
an elliptic pencil $\mid 2F\mid =\mid 2F'\mid$ with $K_S=F-F'$.
 On a general Enriques surface
there are ten different elliptic pencils $\mid 2F_1\mid ,\mid 2F_2\mid,...,
\mid 2F_{10}\mid $ (see [CD]).

\addtocounter{subsection}{1}
\subsection{Mukai lattice}
     It is convenient to describe discrete invariants of sheaves and
bundles on a K3 or an Enriques surface $X$  in the form of vectors in the
algebraic Mukai lattice
  $$M(X) = H^0 (X,\bbbz) \oplus Pic X \oplus H^4 (X,\bbbz) = \bbbz \oplus Pic
X \oplus \bbbz \ni v=(r,D,s)$$
with inner product $< , >$
$$<(r,D,s), (r',D',s')> = rs'+s'r-D.D'$$
To each sheaf $E$ on $X$ with $c_1(E) = D, c_2(E) \in H^4 (X, \bbbz) =
\bbbz$ we associate the vector
$$v(E) = \left( rk(E),D, \frac{1}{2}D^2 - c_2 + rk(E)
\frac{\chi(O_X )}{2}\right) , $$
where $\chi(O_X )=1-q+p_g $ is equal to 2 for a K3 surface and 1 for
an  Enriques surface.
This formula,the Riemann-Roch theorem yield the equalities:
\beq
<v(F),v(E)>=<v(E),v(F)>&=&\chi(E,F)\nonumber\\
 = dim Ext^0 (E,F) -  dim Ext^1 (E,F)&+& dim Ext^2 (E,F)
\enq
 Because $K_X $ is numerically  equal to zero we have
 that $\chi(E,F) = \chi (F,E)$ on a K3 or an Enriques surface.
For the short exact sequence
$$0 \to F \to E \to G \to 0$$
we have the following equalities
$$v(E) = v(F) + v(G) = (r(F)+r(G), c_1 (F) +c_1 (G),s(F) + s(G)).$$

\addtocounter{section}{1}
\subsection{Exceptional bundles}
    {\bf Definition:} E is an exceptional sheaf on surface S if:\\
 ${~~~~~~~~~~~~~~}dim Ext^0 (E,E) = 1, dim Ext^1 (E,E) = q , dim Ext^2 (E,E) =
  p_g$.\newline
     From the Riemann-Roch theorem we have:
$$\chi(E,E) = r^2 \chi(O_S )+(r-1)c_1^2 -2rc_2.$$
Since $H^2 (S,\bbbz )$ is even lattice and $\chi(O_S)=1$ for an Enriques
surface S we see that E has odd rank if E is
an exceptional bundle.
For the description of  exceptional vector bundles we need the following
result of Kuleshov.
\newtheorem{ttt}{Theorem [Ku]}[subsection]
\begin{ttt}
Let X be a smooth a K3 surface  and let  $H$ be an arbitrary ample divisor
on X, and $v=(r,D,s), r > 0$ is an exceptional vector (i.e $v^2 =2$)
belonging to the Mukai lattice on X. Then there exist a simple ,
$\mu_H $-semi-stable bundle E which realize the vector v (i.e.v=v(E)).
\end{ttt}
     I wish to start from some useful facts. The first one is about torsion
free sheaves with some homological condition.
\newtheorem{Mukai}{Proposition}[subsection]
\begin{Mukai}
Let {\it E} be a torsion free sheaf on a smooth  surface S  and $dim Ext^1
({\it E,E}) =1$ or 0. Then {\it E} is locally free.
\end{Mukai}
\pf \     We have the following exact sequence:
$$0 \to {\it E} \to {\it E}^{**} \to M \to 0, $$
where ${\it E}^{**}$ is  double dual of {\it E} and cokernel M is of finite
length. Now I use Mukai [M]  result
 ( see Corollary 2.11 and 2.12) and obtain the following inequality:
$$ dim Ext^1 ({\it E}^{**},{\it E}^{**}) + dim Ext^1
(M,M) \leq  dim Ext^1 ({\it E},{\it E}) $$
Because $v^2 (M)=0$ we have that $ dim Ext^1 (M,M)$ is equal to
 $2dim End_{O_S} (M)$.\\
Since $Ext^1 ({\it E},{\it E})=0$ we
obtain that M=0 and ${\it E} = {\it E}^{**}$. And hence ${\it E}$ is
locally free  the statement follows. $\odot$

     Now I wish formulate the main result.
\newtheorem{exc}[ttt]{Theorem}
\begin{exc}
 Let S be a smooth  Enriques surface S, $v=(r,D,s) \in M(S)\\ ( r > 0 )$ and
 $v^2=1$ then: \\
(i)\ There is an ample divisor H such that $D\cdot H$ and r have not common
divisor greater than 1 (i.e. $(D\cdot H,r) = 1$).   \\
(ii)\ For any ample divisor H with condition  $(D\cdot H,r)=1$ exist an
exceptional
vector bundle E and only one such that v(E)=v and E is H-stable.
\end{exc}
\pf\ \  (i) Because $v^2=1$ we have $2rs-D^2=1$. This means that $D \not\subset
r'\cdot H^2 (S,\bbbz ) $ for any $r'$ such that $( r',r ) > 1 $ (here ( , )
 means the greatest common divisor). Since our
lattice  $H^2 (S,\bbbz ) $ is unimodular there is $X \in H^2 (S,\bbbz ) $
such that $(X\cdot D,r)=1$. Hence $H_k =X + krH$ will be very ample for any
ample divisor H and
$k\gg 0$ . And $H_k$ satisfy our condition $(H_k \cdot D,r)=1$. This prove
the first statement.\\
    (ii) Let a universal covering space  of S be X which is  a K3 surface and
let $\pi$ be the quotient map. Consider the vector $\hat v =(r,D',s')=\pi
^*(v)$ on K3
surface X. It turns out that $\hat v ^2=2$ so by theorem of Kuleshov there is
an exceptional vector bundle F on X such that $\hat v (F)=\hat v $ and F is
$\hat H  = \pi
^*(H)$-semi-stable. In fact F is $\hat H $-stable. Indeed, we have
 $(\hat H \cdot D',r) =(2H\cdot D,r) = 1$ (recall that r is even number
as we notice above) and for
any subsheaf W with $0 < rank(W) < rank(F) $ the following:
 $$\frac{c_1 (W)\cdot
\hat H }{rank(W)}<\frac{D'\cdot \hat H }{rank(F)},$$ therefore F is $\hat H$
-stable .
So F and $\sigma ^* (F)$ both are $\hat H$ -stable, where $\sigma$ is  the free
involution on X such that $X/\sigma = S.$ It is easy to see that
 $\chi(F, \sigma ^* (F)) = 2$. Hence there is non zero homomorphism from F to
$\sigma ^*
(F)$ which should be isomorphism because both vector bundles have the same
determinant and  are $\hat H$ -stable. This isomorphism means that there is a
vector
E on the surface S such that $\pi ^*(E)=F$. Of course, the vector bundle E is
H-stable and v(E)=v.
      Assume that there is another H-stable vector bundle $\hat E$  and
$v(E)=v(\hat E )$. Then $\chi(E,\hat E )=1$, therefore there is non trivial map
from $\phi:E\to \hat E $ or by Serre duality $\rho:\hat E \to E\otimes K$. In
both
cases it should be isomorphisms by stability assumption. But then $E=\hat
E$  or
$det\rho$ gives non zero element of canonical class. This contradiction
prove that H-stable bundle is unique.$\odot$ \\
 {\bf Notice} that from condition
$Ext^1(E,E)=0$ we get only that moduli space of E contain only discrete set of
bundles. But by (ii) such E is only one so moduli space consist of only one
 point.
 \paragraph{Remark}: It is not clear (to  my ) whether an exceptional
H-stable vector bundle  E is G-stable for any another ample divisor G. By
theorem it can happen only if $(G\cdot c_1 (E),rk(E)) > 1$.

\addtocounter{section}{1}
\subsection{Examples}
     I wish to give some explicit examples of \ex\ \vv\ bundles. The one way
to construct them is as in the theorem above to find a stable, invariant,
\ex \vv bundle on a K3 surface.
     Consider  Horikawa's representation of Enriques surface (see for
details [BPV]).
We introduce coordinates $(z_0:z_1:z_2:z_3)$ on ${\bbbp ^3}$ such that a
quadric $Q=\bbbp ^1 \times \bbbp ^1$ is embedded by
$$z_0=x_0y_0,\ \  z_1=x_1y_1,\ \  z_2=x_0y_1,\ \  z_3=x_1y_0.$$
If we define the involution $\tau$ on  ${\bbbp ^3}$ by $\tau
(z_0:z_1:z_2:z_3)=(z_0:z_1:-z_2:-z_3)$, then Q is $\tau$ invariant with
$\tau$ acting on Q by
$$\tau ((x_0:x_1)(y_0:y_1)=(x_0:-x_1)(y_0:-y_1).$$
 On Q the involution $\tau$ has the four fixed points
$$(x_0:x_1)(y_0:y_1)=(1:0)(1:0) , (1:0)(0:1),(0:1)(0:1) , (0:1)(1:0).$$
 Take a polynomial of bidegree (4,4) which define a  $\tau$ invariant curve
 B.   Assume that B have not   any of fixed points. Consider a surface X which
is a double cover of Q
 ramified over B. It turns out  that X is a K3 surface and the involution
 $\tau$ induce the involution $\sigma$ on X which is without fixed points.
 Hence $X/\sigma =S$ is an Enriques surface. Notice that general Enriques
 surface can be obtain by this construction.  Let $\pi :X \to S$ be
 factorization map and $\phi :X \to Q$ the double cover ramified over B.
If E is \ex \vv bundle on Q then $\phi ^* (E)=\hat E $ is an \ex . Indeed , by
projection formula we have $H^i (End_{O_X}(\hat E ) = H^i (End_{O_Q}(E) \oplus
H^i (End_{O_Q}(E) \otimes (-K_Q))$ (recall that $B=-2K_Q$). Since, by the
Serre duality, we have  $H^i (End_{O_Q}(E) \otimes
(-K_Q))=H^{2-i}(End_{O_Q}(E))$ therefore $\hat E$  is an \ex\ \vv\ bundle on
the K3 surface X. This
bundle $\hat E$ is $\sigma$ invariant because E is $\tau$ invariant
 on quadric Q.
(Indeed, E is rigid, therefore E is $PGL(2)\times PGL(2)$-homogeneous.)
Because $\hat E$  is $\sigma$ invariant, there is a \vv \ bundle F on S such
that
$\hat E =\pi ^* (F)$. It is easy to see that $v^2 (F)=1$ and F is rigid,
therefore
F is \ex\ \vv\ bundle on Enriques surface S. I am not able to say anything
about the stability of $\hat E$ and F. This produce a lot of \ex\ \vv
 bundles because we know how to construct all \ex\ \vv\ bundles on smooth
quadric Q.

     In the similar way  we can consider quartic X in $\bbbp ^3$ which is
defined by the equation $z_0^4+z_1^4-z_2^4-z_3^4$. This is a smooth
K3 surface (see for details in [GH]). Let T  be an automorphism on $\bbbp ^3$
defined as follows:
$$T: (z_0,z_1,z_2,z_3) \to (z_0,\sqrt{-1}z_1,-z_2,-\sqrt{-1}z_3).$$
This  automorphism T has 4 fixed points on $\bbbp ^3$ no one of which lay on
the surface X , $T^2$ has 2 fixed lines:
$$l_1=(z_0=z_2=0), \ l_2=(z_1=z_3=0). $$
These lines intersect the surface X in 8 points $p_1,....,p_8.$ Consider
blow-up of X in these 8 points ${\bar X} \to X$. Let ${\bar T}$ denote induced
automorphism on ${\bar X}$. It turns out that
$X'={\bar X}/\{ {\bar T}^{2n} \} $ is a K3 surface and $\bar T$ acts on
$X'$
 as an involution
 without fixed points. So $X'/\bar T$ is an Enriques surface S.
 Now we can get an \ex\
\vv \ bundle on S from any \ex \ \vv bundle on $\bbbp ^3$ because each \ex
\ \vv \ bundle on $\bbbp ^3$ is a homogeneous, therefore it is T invariant.
Of course, we get the \vv bundle on $X'$ which is an \ex  and the bundle
descend to S also as an \ex . This procedure gives us a lot \ex \ \vv \
bundles on S and we can describe it because we know  constructions of \ex
 bundles on $\bbbp ^3$.

    Now I wish  discuss about the ability to construct exceptional
 collections on Enriques surface. Recall that by definition
${E_1,E_2,...,E_n}$ is an  \ex\
 collection if $Ext^i (E_k,E_j)=0$ for any i and $k>j.$ In particular,  we
 have that $\chi (E_k,E_j)=0$. But on an Enriques surface we have $\chi
(E,F)=\chi
 (F,E)$ for any sheaves E and F. Hence for an \ex \ collection on Enriques
 surface should be true the following:
$$ Ext^i (E_k,E_j)=0 , \forall \ i \ and\  k>j;\  \  \ \ \chi (E_a,E_b)=0\ if\
a\not=
b.$$
On a general Enriques surface exist  \ex\ colection with ten bundles. Indeed,
it is well known that
on general Enriques surface there are ten different elliptic pencils say
$\mid 2F_1 \mid ,...,\mid 2F_{10} \mid $ (see [CD]). It is easy to see that
 $Ext^i (F_k,F_j)=0 , \forall \ i \ and\ k\not= j,$
therefore ${F_1,F_2,...,F_{10}}$ is an \ex\ collection.  It will be very
interesting to describe the orthogonal category in the derived category
D(S) (which is finite )
of all sheaves on Enriques surface S. This orthogonal category should have
only two independent elements.

   \addtocounter{section}{1}
\subsection{ Modular operations }
     There are some natural modular operations from one moduli space to
     another which gives  an isomorphism
     of tangent bundles of moduli spaces.For example $E\longleftrightarrow
  E^*, E \longleftrightarrow E\otimes D$, where D is a line bundle.
     On an Enriques surface we have the very interesting modular
operation which I call a reflection. This operation is similar to the
reflection on a K3 surface (see [T] 4.10,4.11). I wish to describe it.
First of all this reflection acts on Mukai lattice in the following way:
$$v=(r,D,s) \longleftrightarrow R(v)=\hat v =\left( 2s,D+\left(
s+\frac{r}{2}\right)
K_S,\frac{r}{2}\right).$$
Notice that $v^2 = 2rs-D^2 = \hat v ^2 $.\\
     Now I describe it on the level of sheaves.
 Assume that a torsion free sheaf E is generated by its section,v(E)=v=(r,D,s)
 and
$\chi(E\otimes K)=h^0 (E\otimes K),\\
 h^1 (E\otimes K)= h^1 (E) = 0$. Notice
that from this we have $h^2 (E\otimes K)=\\
 h^2 (E) = 0, h^0 (E)=\chi
(E)=\chi (E\otimes K)$. For example, if we twist any stable bundle by
sufficiently large  very ample divisor then our condition will be
satisfied.
 Consider the following   exact sequence :
\beg
0 \to \bar E^* \to H^0 (E)\otimes O_S\stackrel{ev}{ \to } E \to 0 ,
\label{ebundle}
\en
where $ev: H^0 (E)\otimes O_S \to E \to 0$ is the canonical evaluation map
(surjective by assumption).
For the convenience denote $H=H^0 (E),\ h=dim H$ and consider the dual
sequence:
\beg
0 \to  E^* \to H^* \otimes O_S \to \bar E \to 0. \label{edual}
\en
By our conditions and Serre duality we have that
$$h^1 (\bar E )=h^2 (E^* )=h^0 (E \otimes K) =\chi (E\otimes K)=h.$$
 Consider the following sequence:
\beg
 0\to H^1(\bar E ) \otimes K_S \to \hat E  \to \bar E \to 0, \label{extention}
\en
where $\hat E$  is given by universal extension element $id \in Ext^1 (\bar E,
H^1(\bar E )) = End(H^1(\bar E )).$ Denote $R(E)=\hat E $.
As an easy consequence of two sequences  (\ref{ebundle}) and
(\ref{extention}) our
assumptions and Serre duality is the
following
\newtheorem{PP}{Proposition}[subsection]
\begin{PP}
Assume E is the sheaf as above then sheaves $\bar E$ and $\hat E $
from sequences (\ref{ebundle}) and (\ref{extention}) satisfies  the following
 properties: \\
1. ${\bar E}$ is globally generated by sections.\\
2. $\chi(\bar E)=0,\ h^0 (\bar E)=h=h^1 (\bar E), \ h^2 (\bar E)=0.$\\
3. $\chi(\bar E\otimes K)=0,\ h^i (\bar E\otimes K_S)=0, \ for \ \forall
i>0$.\\
4. $v(\hat E )=\hat v =R(v)=
\left( 2s,D+\left( s+\frac{r}{2}\right) K_S,\frac{r}{2}\right) \ and\\
 h^0 ( \hat E )= h^0 ( \hat E \otimes K)=h,\ h^i ( \hat E )=
h^i ( \hat E \otimes K)=0 \ for \ \forall i>0$.\\
5.$Hom(E,E)=Hom(\bar E,\bar E)=Hom(\hat E ,\hat E ),Ext^2(\bar E,\bar E)=0.$\\
6. $rank(E)=rank(\hat E)\pmod{2}$. Moreover, if $rank(E)=2k+1$ and
 both E and $\hat E$  are H-stable then
$Ext^2(\hat E ,\hat E )=Ext^2(E,E)=0,Ext^1(\hat E ,\hat E )=Ext^1(E,E)$.
 \label{PP}
\end{PP}
\pf \ \
    1. The sheaf $H^* \otimes O_S$ in the middle of the sequence
    (\ref{edual})
is globally generated by sections so $\bar E$ is too.\\

\noindent 2.The  corresponding long in cohomology to the sequence
(\ref{edual})
 gives us $h^0 (\bar E)=h$ because by Serre duality $h^0
(E^*)=h^2 (E\otimes K)=0$ and $h^1 (E^*)=h^1 (E\otimes K)=0$. Also
$h^2(\bar E)=0$ and $h^1 (\bar E)=h$ as we already noticed. Hence
$\chi(\bar E)=0.$\\

\noindent     3.Consider the sequence (\ref{edual}) twisted by K:
\begin{equation}
0 \to  E^*\otimes K \to H^* \otimes K_S \to \bar E\otimes K \to 0.
\label{etwistK}
\end{equation}
In the same way the corresponding long exact sequence in cohomology gives:
$h^0(\bar E \otimes K)=h^1( E^* \otimes K)=h^1( E) $ and $\ h^1(\bar E \otimes
K)=0$ because $h^2(E^* \otimes K)=h^0(E)=h=h^2(H^*\otimes K);\ h^2(\bar
E\otimes K)
=h^0(\bar E^*)$ (by the sequence (\ref{ebundle})).\\

\noindent     4. An easily calculation shows that
 $v(\hat E )=\hat v =\left( 2s,c_1(E)+\left( s+\frac{r}{2}\right)
K_S,\frac{r}{2}\right).$
 Since the sequence (\ref{extention}) is the
 universal extension we have $h^1( \hat E )=h^2( \hat E )=0$, therefore, by the
 Riemann-Roch theorem,  we obtain that $h^0( \hat E )=h$. If we twist
 (\ref{extention}) by K and
 then use properties of $\bar E\otimes K$ we easily get that\\ $h^0(
 \hat E \otimes K)=h,\ h^1(\hat E \otimes K)=h^2(\hat E \otimes K)=0$. \\

\noindent     5. Applying $Hom(\bar E,*)$ to (\ref{edual}) we get the long
exact sequence:
\begin{eqnarray*}
0\to Ext^0(\bar E,E^*)\to & Ext^0(\bar E,H^*\otimes O_S) & \to Ext^0(\bar
E,\bar E) \to \\
\to Ext^1(\bar E,E^*)\to & Ext^1(\bar E,H^*\otimes O_S) & \to Ext^1(\bar
E,\bar E) \to \\
\to Ext^2(\bar E,E^*)\to & Ext^2(\bar E,H^*\otimes O_S) & \to Ext^2(\bar
E,\bar E) \to0
\end{eqnarray*}
By Serre duality and the statement 3 Ext groups in the middle are\\
 $H^i(\bar E\otimes K)\otimes
H^* =0$. Hence we have $ Ext^1(\bar E,E^*)= Ext^0(\bar E,\bar E) ,\\
Ext^2(\bar E,\bar E) =0.$
Also applying $Hom(*,E^*) $ to (\ref{edual}) we get:
\begin{eqnarray*}
0\to Ext^0(\bar E,E^*)\to & Ext^0(H^*\otimes O_S,E^*) & \to Ext^0(E^*,E^*) \to
\\
\to Ext^1(\bar E,E^*)\to & Ext^1(H^*\otimes O_S,E^*) & \to Ext^1(E^*,E^*)
\to
\end{eqnarray*}
By our assumption and Serre duality the middle Ext groups are \\
$Ext^0(H^*\otimes O_S,E^*)=H^2(E\otimes K) = 0 = Ext^1(H^*\otimes
O_S,E^*)=H^1(E\otimes K)$.\\
 Hence $Hom(E^*,E^*)=Hom(E,E)=Hom(\bar E,\bar E).$
Now applying $Hom(\hat E ,*)$ to (\ref{extention}) we get:
\begin{eqnarray*}
0\to Ext^0(\hat E ,H\otimes K)\to & Ext^0(\hat E ,\hat E ) & \to Ext^0(\hat E
,\bar E) \to \\
0\to Ext^1(\hat E ,H\otimes K)\to & Ext^1(\hat E ,\hat E ) & \to Ext^1(\hat E
,\bar E) \to
\end{eqnarray*}
Since $Ext^0(\hat E ,H\otimes K)=H^2(\hat E )\otimes H=0$ we obtain that
$Ext^0(\hat E ,\hat E )= Ext^0(\hat E ,\bar E)$. And, in the similar way, after
applying
$Hom(*,\bar E)$, we get that $Ext^0(\hat E ,\hat E )= Ext^0(\bar E,\bar
E)$.

     6.If rank(E) is odd then $s(E)\in \frac{1}{2}\bbbz$ but
     $s(E)\not\subset \bbbz $, therefore $rank(\hat E)\\
 = 2s(E)$
 is odd number too. If rank(E) is even then $s(E)\in \bbbz$, therefore
 $rank(\hat E) = 2s(E)$
 is even.Consider any non zero element $\phi$ in $Ext^2(E,E)=\\Ext^0(E,E\otimes
K)^*$
by the stability assumption this should be an isomorphism but then $det\phi$
is non zero element of canonical class. This contradiction shows that
$Ext^2(E,E)=0$ and $Ext^2(\hat E ,\hat E )=0$. Because $v^2=\hat v ^2$ we have
$Ext^1(E,E)=Ext^1(\hat E ,\hat E )$.$\odot$\\

     I am able to reverse this operation in the following situation.
     Consider   a sheaf F and the following exact sequence:
\beg
0\to H^0( F\otimes K_S ) \otimes K_S \stackrel{ev}{\to} F \to \bar F \to 0,
 \label{fx}
\en
where $ ev:H^0( F\otimes K_S ) \otimes K_S \to F$ is the canonical  evaluation
map.
Assume that $\bar F$ is globally generated a torsion free sheaf then we have
the
following exact sequence:
\beg
0 \to  \hat F ^* \to H^0 (\bar F) \otimes O_S \stackrel{ev}{\to} \bar
F \to 0. \label{ff}
\en
Denote $R(F)=\hat F$.
Under these assumptions we can prove the similar result:
\newtheorem{P'}[PP]{Proposition}
\begin{P'}
Assume F is the sheaf as above then sheaves $\bar F$ and $\hat F$  from
sequences
(\ref{ff}) and (\ref{fx}) satisfies  the following  properties: \\
1. $\chi(\bar F)=0,\ h^0 (\bar F)=h=h^1 (\bar F), \ h^2 (\bar F)=0.$\\
2. $\chi(\bar F\otimes K)=0,\ h^i (\bar F\otimes K_S)=0, \ for \ \forall
i>0$.\\
3. $v(\hat F )=\hat v =\left( 2s,D+\left( s+\frac{r}{2}\right)
 K_S,\frac{r}{2}\right) \ and\  h^0 ( \hat F )= h^0 ( \hat F \otimes K)=h,\\
 \ \ \ h^i ( \hat F )=
{~~}h^i ( \hat F \otimes K)=0 \ for \ \forall i>0$.\\
4.$Hom(F,F)=Hom(\bar F,\bar F)=Hom(\hat F ,\hat F ),Ext^2(\bar F,\bar F)=0.$\\
5. $rank(F)=rank(\hat F)\pmod{2}$. Moreover, if $rank(F)=2k+1$ and
 both F and $\hat F$  are H-stable then
$Ext^2(\hat F ,\hat F )=Ext^2(F,F)=0,Ext^1(\hat F ,\hat F )=Ext^1(F,F)$.
 \label{P'}
\end{P'}
Notice that from both propositions follows that R(R(E))=E.
\paragraph{Remark:} The reflection always exist in the derived category of
sheaves
on surface S. It does not matter whether the evaluation map
$  H^0( F ) \otimes O_S \stackrel{ev}{\to} F$ is a surjective or an
injective map.\\

     Now I wish to give a few examples of reflections.\\
1. Assume we have a smooth curve C on S,  A a globally generated
divisor on the curve C with the properties $h^1 (O_C(A))=h^1(O_C(A\otimes
K_S))=0$.
We consider the following exact sequence:
\begin{equation} 
0 \to E(C,A)^* \to H^{0}(O_C(A)) \otimes O_S\stackrel{ev}{\to} O_{C}(A)
 \to 0  \label{E*}
\end{equation}
The dual sequence to (\ref{E*}) is :
\begin{equation}
0 \to H^{0}(A)^{*} \otimes O_S \to E(C,A) \to O_{C}(C) \otimes  A^* \to 0
\end{equation}
By our assumption, this sequence and Serre duality on the curve C,
we obtain that
$h^1(E(C,A))=h^1(O_C(C-A))=h^0(O_C(A)$ and get
the following sequence:
\begin{equation} 
0 \to H^1(E(C,A))\otimes K_S \to \hat E  \to E(C,A)  \to 0
\end{equation}
So we have $R(O_C(A))=\hat E $. \\
If C is (-2)-curve and $A=O_C$,
we  see that $R(O_C)=F$ is an extremal rank 2 vector bundle (i.e.
$Ext^0(F,F)=Ext^2(F,F)=\bbbc, Ext^1(F,F)=0.)$

2.If $\mid 2F\mid,\ \mid 2G\mid$ are two elliptic pencils on S then a  pencil
$\mid F+G\mid$ has two different base points x and y. From the standard
sequence:
\begin{equation}
0 \to  J_{x+y} (F+G) \to O(F+G)\to O_{x+y}(F+G)\to  0, \label{jj}
\end{equation}
 we see that $h^1(J_{x+y} (F+G))=2$. This gives us the following bundle E,
 defined by the universal extension element $id\in End(H^1(J_{x+y} (F+G))$:
\begin{equation}
0 \to H^1(J_{x+y} (F+G))  \otimes K_S \to E\to J_{x+y} (F+G) \to 0, \label{ej}
\end{equation}
By lemmas 1.1,1.2 in [T], E is a simple bundle. An easily calculation
shows that E is an \ex \ bundle. Also by \ref{P'}, we see that $R(E)=O(F+G)$.
Notice that $O(F+G)$ is not globally generated by section, therefore we cannot
use \ref{PP} to this bundle to produce $R(O(F+G))$, but we can  do this in
the derived category. Since $R(E)=O(F+G)$, we also obtain  that
$R(O(F+G))=E$.

3.If we consider a divisor $O(aF+bG)$ for $a\geq b\geq 2$ then this divisor
will be an ample on a general Enriques surface. Hence $R(O(aF+bG))$ is an \ex\
bundle of rank $2ab+1$.
\paragraph{References} : \newline

 [BPV] W.Bart, C.Peter, A.Van de Ven. Compact complex surfaces.
 Berlin, Heidelberg,New York :Springer 1984. \newline

[CD] F. Cossec ; I. Dolgachev "Enriques Surfaces 1", Birkh\"{a}user 1989.
 \newline

[GH] Griffiths,Ph., Harris J.,Principles of algebraic geometry, New
York\newline
 (1978) \newline

[Ki1] Kim,Hoil.:Exceptional bundles on nodal Enriques surfaces,\newline
 Bayreuth
preprint(1991). \newline

[Ki2] Kim,Hoil.:Exceptional bundles and moduli spaces of stable vector
\newline
 bundles on  Enriques surfaces,Bayreuth
preprint(1991). \newline

[Ku] Kuleshov,S.A:An existence theorem for \ex\ bundles on K3 \newline
 surfaces,Math.USSR Izvestia,vol.34,373-388.(1990). \newline

[M] Mukai,S.On the moduli space of bundles on K3 surfaces I, in Vector
\newline
 Bundles, ed. Atiay et all, Oxford University Press, Bombay, 341-413(1986).
\newline

[N] Daniel Naie. Special rank two vector bundles over Enriques surfaces,
\newline
preprint. \newline

[T] A.N. Tyurin "Cycles, curves  and vector bundles on algebraic surfaces."
Duke Math.J. 54, 1-26,(1987).  \newline

Department of geometry and topology,

Faculty of mathematics,

Vilnius university,

Naugarduko g.24,

2009 Vilnius, Lithuania.

e-mail:Severinas.Zube@maf.vu.lt
\end{document}